\documentclass[11pt]{article}
\usepackage{amsmath}
\usepackage{amsfonts}
\usepackage{amssymb}
\usepackage{amsthm}

\newcommand{\Z}{{\bf{Z}}}
\newcommand{\Q}{{\bf{Q}}}

\newcommand{\R}{{\bf{R}}}
\newcommand{\C}{{\bf{C}}}

\newcommand{\bolda}{{\bf{a}}}

\newcommand{\ma}{{\mathfrak{a}}}

\newcommand{\ra}{\rightarrow}

\newcommand{\aaa}{{\mathfrak{a}}}

\newcommand{\inti}{{\int_0^\infty}}

\newtheorem{lem}{Lemma}[section]

\newtheorem{thm}[lem]{Theorem}

\theoremstyle{definition}

\newcommand{\later}[1]{}
\newcommand{\comment}[1]{}
\newcommand{\com}[1]{}

\newcommand{\thetitle}
{A generalization of Kronecker's first limit formula to GL(n)}

\begin{document}
\parindent=2em

\date{}

\title{\thetitle}
\author{Amod Agashe\footnote{The author was partially supported by the 
Collaboration Grants for Mathematicians from the Simons Foundation.}}
\maketitle


\begin{abstract}
Kronecker's first limit formula gives 
the polar and constant terms of the Laurent
series expansion of the Eisenstein series for~${\rm SL}(2,\Z)$
at~$s=1$, which in turn can be used to find expressions for
the polar and constant terms of 
the Dedekind zeta function of quadratic fields.
In this article, we generalize the formula to certain maximal parabolic
Eisenstein series associated to~${\rm SL}(n,\Z)$ for~$n \geq 2$.
We also show how the generalized
formula can be used to find expressions for the polar and constant
terms of the Dedekind zeta function of any number field at~$s=1$.
\end{abstract}

\section{Introduction}
Let $n \geq  2$ be an integer.
Let $\tau$ be in the generalized upper half-plane~$\mathfrak{H}^n$,
which consists of $n \times n$ matrices with real number entries
that are the product of an 
upper triangular matrix with $1$'s along the diagonal and a diagonal matrix
with positive diagonal entries such that the lowermost 
diagonal entry is~$1$. 
When $n=2$, one can identify~$\mathfrak{H}^n$ with the usual complex
upper half plane (for details, see, e.g.,~\cite[\S1.2]{goldfeld}).

In the following, $m_1, \ldots, m_n$ denote integers with 
perhaps some added restrictions as noted;
in particular, we follow the convention
that in any sum over a subset of~$m_1, \ldots, m_n$, 
if a term has denominator zero
for some values of $m_1, \ldots, m_n$,
then the term is to be skipped in the sum. 

Consider the maximal parabolic Eisenstein series
$$E_n(\tau, s)  = \sum_{(m_1,\ldots, m_n) = 1} \frac{(\det \tau)^s}
{\parallel (m_1 \ldots m_n) \tau \parallel^{ns/2}},$$
where ${\parallel (m_1 \ldots m_n) \tau \parallel}$ denotes
the norm of the row vector that is the product of the row vector
$(m_1 \ldots m_n)$ and the matrix~$\tau$; this series converges 
when ${\rm Re}(s) > 1$, and is known to have a 
meromorphic continuation to all of~$\C$.

Let
\begin{eqnarray} \label{eqn:e*}
E_n^*(\tau, s) 
& = &  
\pi^{-ns/2} \Gamma(ns/2) \zeta(ns) E(\tau, s) \nonumber \\
& = & \pi^{-ns/2} \Gamma(ns/2) 
{\sum_{m_1,\ldots, m_n}} 
\frac{(\det \tau)^s}{\parallel (m_1 \ldots m_n) \tau \parallel^{ns/2}} \ .
\end{eqnarray}

Note that in the case where $n=2$, 
if $\tau$ corresponds to the point $z = x + iy$ in the complex upper half
plane, then
$$
E^*_2(\tau, s) 
= \pi^{-s} \Gamma(s) {\sum_{m_1, m_2}} 
\frac{y^s}{|m_1 z  +  m_2|^{2s}}\ . 
$$
The classical Kronecker's first limit formula 
gives the first two terms of the Laurent expansion of~$E^*_2(\tau,s)$
at $s=1$: 
\begin{eqnarray} \label{eqn:kron}
E^*_2(\tau, s)
= \frac{1}{s-1} + \big( \gamma  - \log 4 \pi
- \log y
- 4 \log |\eta(\tau)|\big) + O(s-1) \ ,
\end{eqnarray}
where $\gamma$ is the Euler-Mascheroni constant
and $\eta(z)$ is the Dedekind eta-function.

Recall that if $K$ is a number field and 
$A$ is an ideal class of~$K$, then the partial zeta function
associated to~$A$ is given by
$$\zeta_K(s,A) = \sum_{\aaa \in A} \frac{1}{{N\aaa}^s}.$$
The Dedekind zeta function is the sum of the partial
zeta functions over all ideal classes. For a number field~$K$, 
we denote by~$w_K$ the number of roots of unity in~$K$ 
and by~$d_K$ the discriminant of~$K$. 

Using~(\ref{eqn:kron}), 
Kronecker showed that if $K$ is a quadratic imaginary
field, then
$$ \zeta_K(s,A)  = \frac{1}{w_K} \frac{2 \pi}{\sqrt{d_K}}
\bigg( \frac{1}{s-1} + 2 \gamma 
- \log 2 - \log y  - 4 \log(|\eta(\tau)|) \bigg) + O(s-1) \ ,
$$
where $\tau$ is an element of the upper half plane such
that $\{1, \tau\}$ is a basis for an ideal in the inverse class of~$A$
and $y$ is the imaginary part of~$\tau$.

Later, Hecke used Kronecker's limit formula~(\ref{eqn:kron})
to show that if $K$ is a real quadratic field and $\epsilon$
is a fundament unit of~$K$, then
\begin{eqnarray} \label{eqn:hecke}
\frac{1}{2}  \big(\pi^{-1} d_K^{1/2}\big)^s \Gamma(s/2)^2 
\zeta_K(s,A)  & =  &
\frac{\log \epsilon}{s-1}
+ (\gamma - \log 4 \pi) \log \epsilon  
- \int_1^\epsilon
\log y(t) \frac{dt}{t} \nonumber
\\ & & 
- 4 \int_1^\epsilon \log |\eta(\tau(t))|
\frac{dt}{t}  + O(s-1) \ ,
\end{eqnarray}
where $\tau(t)$ is a certain 
curve in the upper half plane (which depends on~$A$),
and $y(t)$ denotes its $y$-coordinate.

Kronecker's first limit formula~(\ref{eqn:kron})
and Hecke's formula~(\ref{eqn:hecke}) 
were generalized to the case $n=3$ in~\cite{efrat} 
(see also~\cite{bump-goldfeld} for
a slightly different generalization in this case).
In this article, we generalize 
Kronecker's first limit formula~(\ref{eqn:kron})
and Hecke's formula~(\ref{eqn:hecke}) to arbitrary $n \geq 2$
(see Theorem~\ref{thm:main} and Theorem~\ref{thm:zeta}). 
A generalization of Kronecker's first limit formula 
to arbitrary $n \geq 2$
is also given in~\cite{terras} for Epstein zeta functions.

\begin{thm} \label{thm:main}
For a given $\tau \in \mathfrak{H}^n$,
let $y_1, \ldots, y_{n-1}$ denote the 
unique positive real numbers  such that
for $i \geq 1$, we have $\tau_{n-i,n-i} = \prod_{j=1}^{i} y_j$.
If $m_1, \ldots, m_n$ are integers, then
for $j = 2, \ldots, n$, let 
$b_{j} = \sum_{i = 1, \ldots, j-1} m_i \tau_{i,j}$ and
$c_{j} = \tau_{j,j}$; also let 
$m$ be the nonnegative real number such that
$m^2 = m_n^2/c_n^2 + m_{n-1}^2/c_{n-1}^2 + \cdots + m_2^2/c_2^2$
and 
$d = b_n m_n/c_n + b_{n-1} m_{n-1}/c_{n-1} + \cdots b_2 m_2/c_2$. 
Let $\tau'$ be the submatrix of~$\tau$ obtained by
removing the topmost row and leftmost column and let
$$g(\tau) = \exp \Bigg\{ -\frac{1}{4}  \Bigg( 
\bigg( \prod_{i=1}^{n-1} y_i^{(\frac{i}{n-1})} \bigg)
E^*_{n-1}\bigg(\tau', \frac{n}{n-1}\bigg)  + 
\sum_{m_1 \neq 0} \frac{1}{|m_1|}\sum_{\stackrel{(m_2, \ldots,  m_{n})}{\neq (0,\ldots, 0)}} 
 \exp\Big(2 \pi i d- 2 \pi  |m_1| m \prod_{i=1}^{n-1} y_i\Big) \Bigg) \Bigg\} .
$$
Then
\begin{eqnarray} \label{eqn:kron1}
E_n^*(\tau, s)  & =  &
\frac{2/n}{s-1} + \bigg(\gamma - \log 4\pi 
- \frac{2}{n}\log \bigg(\prod_{i=1}^{n-1} y_i^{i}\bigg) - 4 \log g(\tau)
\bigg)  + O(s-1) \ .
\end{eqnarray}
\end{thm}

The theorem is proved in Section~\ref{sec:proof}.
Our proof is a generalization of the proof of the classical 
Kronecker's first limit formula given in~\cite[\S20.4]{lang:elliptic}, and
the key observation is to break the sum in~(\ref{eqn:e*})
over $m_1, \ldots, m_n$ into two parts conveniently (we break it 
as a sum over~$\{m_1\}$ and a sum over~$\{m_2, \ldots, m_n\}$;
when $n=2$, there is not much of a choice)
and to apply the Poisson summation formula in the 
second sum in the correct order (over 
$m_n$ first, followed by~$m_{n-1}$, and so on, up to~$m_2$; 
i.e., in reverse order).
Even in the case $n=3$, our proof differs in some key steps
with that in~\cite{efrat} (who introduces complex coordinates on~$\mathfrak{H}^3$,
while we don't) and that in~\cite{bump-goldfeld} (who use minimal parabolic
Eisenstein series). Our proof techniques are similar to those used 
in~\cite{terras} (about which we learned only after a first draft of this
article was written), but are
more direct and elementary (the main goal
of~\cite{terras} is to prove the functional equation of the Epstein zeta function
using generalizations of the Selberg-Chowla formula).

If we put $n=2$ in the formula for~$E^*_n(\tau, s)$ above,
then we get the classical Kronecker's first limit
formula~(\ref{eqn:kron}), with $g(\tau) = |\eta(\tau)|$
(this last equality follows from the correctness of our formula,
but also because our proof is identical to the classical proof
in the case $n=2$ as given in~\cite[\S20.4]{lang:elliptic}).
Thus our function~$g(\tau)$ for arbitrary~$n$
is a generalization of $|\eta(\tau)|$.
If we put $n=3$, then we recover the formula for $E_3^*(\tau, s)$
given in~\cite[Theorem~1]{efrat}.

\comment{
We suspect that the expression for~$E_n^*(\tau, s)$ above can be 
used to show that~$g(\tau)$ is automorphic and
that  $\log (g(\tau))$ is a harmonic function on~$\mathfrak{H}^n$
(i.e., is annihilated by invariant differential operators),
as is done for $n=2$ (see, e.g., \cite[\S~I.2]{siegel:adv} for automorphy)
and~$n=3$ (see~\cite[\S3]{efrat}). This may be the approach to answer
the question raised at the end of $\S1$ and~$\S3$ in~\cite{terras}.
However,
we shall not pursue these issues in the present article.
}

The classical Kronecker limit formula has several applications (see, e.g., \cite{siegel:adv}); many of these generalize to give applications
of our generalization of the limit formula. In this article, we shall
limit ourselves to one application:
in Section~\ref{sec:zeta}, 
we a formula analogous to~(\ref{eqn:hecke}) 
for the polar and constant terms of the Laurent series
of the partial or Dedekind zeta function
for any number field: see
Theorem~\ref{thm:zeta}
(as mentioned earlier, 
this was already done for cubic fields in~\cite{efrat}).

\section{Zeta functions of number fields} \label{sec:zeta}

In this section, we  a formula for 
the Laurent series
expansion of the partial and Dedekind zeta function of a number field
times some explicit functions; 
the main idea is to use a generalization of a trick of Hecke
to express the partial zeta function as an integral of an Eisenstein series
over a suitable region and then use 
our limit formula for the Eisenstein series. 
The procedure is described for $n=3$ in~\cite[\S4]{efrat}, 
and the discussion below is its generalization.
However, we do need a new input, which is formula~(\ref{eqn:complex}),
the analog of which was not required in loc. cit.

Let $K$ be a number field of degree~$n$ over~$\Q$. 
Let $r$ denote the number of real embeddings and $c$ denote the number
of complex conjugate embeddings. We assume that $r+c > 1$ since the cases 
where $r+c=1$ are classical (the field of
rational numbers and quadratic imaginary
fields). 

Let $A$ be an ideal class of~$K$. 
Fix $B \in A^{-1}$. Let  $U$ denote the unit group of~$K$
and as before, let $w_K$ denote the number of roots of unity in~$K$.
Let $m = r + c -1$, and let $\epsilon_1, \ldots, \epsilon_m$
denote a fundamental set of units. 
Then 
\begin{eqnarray}\label{eqn:zeta}
\zeta_K(s,A) = \sum_{\aaa \in A} \frac{1}{{N\aaa}^s} = 
{NB}^s \sum_{\lambda \in B/U} \frac{1}{|N\lambda|^s} =
\frac{{NB}^s}{w_K} \sum_{\lambda \in B/\langle \epsilon_1, \ldots, \epsilon_m \rangle} \frac{1}{|N\lambda|^s} .
\end{eqnarray}

Order the embeddings of~$K$ so that the first $c$ are complex and the remaining~$r$ are real. 
If $x \in K$, then 
for $i = 1, \ldots, r+c$, let $|x|_i$ denote the absolute value of
the image of~$x$ under the $i$-th embedding. 
For $i = 1, \ldots, m+1 = r+c $, 
let $\delta_i = 1$ if
the $i$-th embedding is real and $\delta_i =2$ otherwise. 

We first deal with the case where $K$ has at least one real embedding.
Then the $(m+1)$-st embedding is real, and thus $\delta_{m+1} = 1$. 
A change of variables shows that
for any positive real numbers $a_1, \ldots, a_{m+1}$, 
\begin{eqnarray} \label{eqn:real}
\lefteqn{\int_0^\infty \cdots \int_0^\infty (a_1^2 t_1^2 + \ldots + a_m^2 t_m^2
+ a_{m+1}^2 (t_1^{\delta_1} \cdots t_m^{\delta_m})^{-2})^{-ns/2} 
\frac{dt_1}{t_1} \cdots \frac{dt_m}{t_m}} \nonumber \\
& & = (a_1^{\delta_1} \cdots a_{m+1}^{\delta_{m+1}})^{-s}
\int_0^\infty \cdots \int_0^\infty (t_1^2 + \ldots + t_m^2
+ (t_1^{\delta_1} \cdots t_m^{\delta_m})^{-2})^{-ns/2} 
\frac{dt_1}{t_1} \cdots \frac{dt_m}{t_m}  \ .
\end{eqnarray}
Now for $i = 1, \ldots, r+c$, let $a_i=|\lambda|_i$. 
Then $a_1^{\delta_1} \cdots a_{m+1}^{\delta_{m+1}} = |N\lambda|$.
Let 
\begin{eqnarray} \label{eqn:d1}
d(s) = \int_0^\infty \cdots \int_0^\infty (t_1^2 + \ldots + t_m^2
+ (t_1^{\delta_1} \cdots t_m^{\delta_m})^{-2})^{-ns/2} 
\frac{dt_1}{t_1} \cdots \frac{dt_m}{t_m}  
= \frac{\Gamma(\delta_1 s/2) \cdots \Gamma(\delta_{m+1} s/2)}{2^m (1+ \delta_1
+ \cdots + \delta_m) \Gamma(ns/2)}
\ .
\end{eqnarray}
Putting all this in formula~(\ref{eqn:real}), we see that 
$$
d(s) \frac{1}{|N\lambda|^s} = 
\int_0^\infty \cdots \int_0^\infty (|\lambda|_1^2 t_1^2 + \ldots + 
|\lambda|_m^2 t_m^2
+ |\lambda|_{m+1}^2 (t_1^{\delta_1} \cdots t_m^{\delta_m})^{-2})^{-ns/2} 
\frac{dt_1}{t_1} \cdots \frac{dt_m}{t_m}, $$ 
and so
\begin{eqnarray} \label{eqn1}
& & d(s) \sum_{\lambda \in B/\langle \epsilon_1, \ldots, \epsilon_m \rangle} 
\frac{1}{|N\lambda|^s} \nonumber \\
& = &
\sum_{\lambda \in B/\langle \epsilon_1, \ldots, \epsilon_m \rangle}
\int_0^\infty \cdots \int_0^\infty (|\lambda|_1^2 t_1^2 + \ldots + 
|\lambda|_m^2 t_m^2
+ |\lambda|_{m+1}^2 (t_1^{\delta_1} \cdots t_m^{\delta_m})^{-2})^{-ns/2} 
\frac{dt_1}{t_1} \cdots \frac{dt_m}{t_m} 
\end{eqnarray}
Now $1 = |N \epsilon_i| = \prod_{j=1}^{m+1} |\epsilon_i|_j^{\delta_j}$,
and so 
\begin{eqnarray} \label{eqn:units}
|\epsilon_i|_{m+1} = (|\epsilon_i|_1^{\delta_1} \cdots
|\epsilon_i|_m^{\delta_m})^{-1}, 
\end{eqnarray}
considering that $\delta_{m+1} = 1$. 
For $i = 1, \ldots, m$, let $\epsilon_i$ act on~$(\R^\times_+)^m$ 
by multiplying the $j$-th coordinate by~$|\epsilon_i|_j$.
Letting $D$ be a fundamental domain
under the action of~$\langle \epsilon_1, \ldots, \epsilon_m \rangle$ 
on~$(\R^\times_+)^m$ , we get
\begin{eqnarray} \label{eqn:2}
\lefteqn{
\sum_{\lambda \in B/\langle \epsilon_1, \ldots, \epsilon_m \rangle}
\int_0^\infty \cdots \int_0^\infty (|\lambda|_1^2 t_1^2 + \ldots + 
|\lambda|_m^2 t_m^2
+ |\lambda|_{m+1}^2 (t_1^{\delta_1} \cdots t_m^{\delta_m})^{-2})^{-ns/2} 
\frac{dt_1}{t_1} \cdots \frac{dt_m}{t_m}} \nonumber \\
& = &
\sum_{\lambda \in B}
\int_D  (|\lambda|_1^2 t_1^2 + \ldots + |\lambda|_m^2 t_m^2
+ |\lambda|_{m+1}^2 (t_1^{\delta_1} \cdots t_m^{\delta_m})^{-2})^{-ns/2} 
\frac{dt_1}{t_1} \cdots \frac{dt_m}{t_m} 
\end{eqnarray}

Let $\alpha_1, \ldots, \alpha_n \in K$ be a $\Z$-basis of~$B$.
Let $\lambda \in B$. Then $\lambda = m_1 \alpha_1 + \cdots
m_n \alpha_n$ for some $m_1, \ldots, m_n \in \Z$.
For $i = 1, \ldots, r+c$, let $\lambda_i$ denote the image of~$\lambda$ 
under the $i$-th embedding. We denote $(\alpha_j)_i$ by~$\alpha_{j,i}$.
Recalling the way we ordered the embeddings of~$K$, we see that
for $i = 1, \ldots, c$,
$|\lambda|_i^2 = \Re(\lambda_i)^2 + \Im(\lambda_i)^2$, while
for $i = c+1, \ldots, c+r$,
$|\lambda|_i^2 = \lambda_i^2$
(if $c=0$, then any expression below containing $\Re$ or~$\Im$ 
should be ignored). Thus
\begin{eqnarray*}
& & 
|\lambda|_1^2 t_1^2 + \ldots + |\lambda|_m^2 t_m^2
+ |\lambda|_{m+1}^2 (t_1^{\delta_1} \cdots t_m^{\delta_m})^{-2}\\
& = & 
(m_1 \Re(\alpha_{1,1}) + \ldots + m_n \Re(\alpha_{n,1}))^2 t_1^2 +
(m_1 \Im(\alpha_{1,1}) + \ldots + m_n \Im(\alpha_{n,1}))^2 t_1^2 +
\ldots  \\
&+ & 
(m_1 \Re(\alpha_{1,c}) + \ldots + m_n \Re(\alpha_{n,c}))^2 t_1^2 +
(m_1 \Im(\alpha_{1,c}) + \ldots + m_n \Im(\alpha_{n,c}))^2 t_c^2 \\
&+ & (m_1 \alpha_{1, c+1} + \ldots + m_n \alpha_{n, c+1})^2 t_{c+1}^2 + 
\ldots \\
&+ & (m_1 \alpha_{1, m} + \ldots + m_n \alpha_{n, m})^2 t_{m}^2 \\
&+ & (m_1 \alpha_{1, m+1} + \ldots + m_n \alpha_{n, m+1})^2 
(t_1^{\delta_1} \cdots t_m^{\delta_m})^{-2}\\
& = & 
{\bf m} M M^T {\bf m}^T, 
\end{eqnarray*}
where ${\bf m}$ is the row vector with entries
$m_1, \ldots, m_n$ and 
$M$ is the $n \times n$ matrix whose 
$i$-th row has entries
$$\Re(\alpha_{i,1})t_1, \ \Im(\alpha_{i,1})t_1, \ 
\ldots, \ \Re(\alpha_{i,c})t_c, \    \Im(\alpha_{i,c})t_c, \  
\alpha_{i,c+1} t_{c+1}, \ \ldots, \ 
\alpha_{i,m} t_{m},  \  
\alpha_{i,m+1} (t_1^{\delta_1} \cdots t_m^{\delta_m})^{-2}\ .
$$ 
Let $Q = M M^T$. Then $Q$ is an $n \times n$ positive definite
symmetric matrix, and thus can be written as
$Q = (\det Q)^{1/n} (\det \tau(t_1, \ldots, t_m))^{-2/n} \tau(t_1, \ldots, t_m) 
\tau(t_1, \ldots, t_m)^T$
for some uniquely defined $\tau: \R^m \ra \mathfrak{H}^n$ (note that
$\det Q$ is independent of~$t_1, \ldots, t_m$). 
Thus 
\begin{eqnarray*}
& & 
\sum_{\lambda \in B}
\int_D  (|\lambda|_1^2 t_1^2 + \ldots + |\lambda|_m^2 t_m^2
+ |\lambda|_{m+1}^2 (t_1^{\delta_1} \cdots t_m^{\delta_m})^{-2})^{-ns/2} 
\frac{dt_1}{t_1} \cdots \frac{dt_m}{t_m} \\
& = & 
\sum_{m_1, \ldots, m_n} \int_D ({\bf m} M M^T {\bf m}^T)^{-ns/2}
\frac{dt_1}{t_1} \cdots \frac{dt_m}{t_m} \\
& = & 
\sum_{m_1, \ldots, m_n} \int_D ({\bf m} Q {\bf m}^T)^{-ns/2}
\frac{dt_1}{t_1} \cdots \frac{dt_m}{t_m} \\
& = & 
\sum_{m_1, \ldots, m_n} \int_D ({\bf m} 
(\det Q)^{1/n} (\det \tau(t_1, \ldots, t_m) )^{-2/n} \tau(t_1, \ldots, t_m) 
\tau(t_1, \ldots, t_m)^T {\bf m}^T)^{-ns/2}
\frac{dt_1}{t_1} \cdots \frac{dt_m}{t_m} \\
& = & 
(\det Q)^{-s/2} 
\int_D 
(\det \tau(t_1, \ldots, t_m) )^{s}
\sum_{m_1, \ldots, m_n}  ({\bf m} 
 \tau(t_1, \ldots, t_m) 
\tau(t_1, \ldots, t_m)^T {\bf m}^T)^{-ns/2}
\frac{dt_1}{t_1} \cdots \frac{dt_m}{t_m} \\
& = & 
(\det Q)^{-s/2} \pi^{ns/2} \Gamma(ns/2)^{-1}
\int_D E^*_n(\tau(t_1, \ldots, t_m), s) 
\frac{dt_1}{t_1} \cdots \frac{dt_m}{t_m},
\end{eqnarray*}
by definition of~$E^*_n(\tau, s)$ (see formula~(\ref{eqn:e*})).
Thus from the equation above and equations~(\ref{eqn:zeta}), (\ref{eqn1}), and~(\ref{eqn:2}),
we see that
$$d(s) \zeta_K(s,A) 
= \frac{{NB}^s}{w_K} 
(\det Q)^{-s/2} \pi^{ns/2} \Gamma(ns/2)^{-1}
\int_D E^*_n(\tau(t_1, \ldots, t_m), s) 
\frac{dt_1}{t_1} \cdots \frac{dt_m}{t_m}.$$
Considering that for a complex number~$z$,
$\Re(z) = (z + \overline{z})/2$ and 
$\Im(z) = (z - \overline{z})/2$, we see that 
$4^{c} \det Q = {\rm disc} B = (NB)^2 d_K$, where recall that 
$d_K$ denotes the discriminant of~$K$. So
\begin{eqnarray} \label{eqn:seminfinal}
d(s) \zeta_K(s,A) 
= \frac{(2^{c} d_K^{-1/2})^s }{w_K} \pi^{ns/2} \Gamma(ns/2)^{-1}
\int_D E^*_n(\tau(t_1, \ldots, t_m), s) 
\frac{dt_1}{t_1} \cdots \frac{dt_m}{t_m}.
\end{eqnarray}

Now consider the case where $K$ does not have a real embedding, i.e.,
$K$ is totally complex. In that case, formula~(\ref{eqn:real})
does not work, and instead, we use the following formula:
\begin{eqnarray} \label{eqn:complex}
\lefteqn{\int_0^\infty \cdots \int_0^\infty (a_1^2 t_1^2 + \ldots + a_m^2 t_m^2
+ a_{m+1}^2 (t_1 \cdots t_m)^{-2})^{-ns/2} 
\frac{dt_1}{t_1} \cdots \frac{dt_m}{t_m}}  \nonumber \\
& & = (a_1^2 \cdots a_{m+1}^2)^{-s}
\int_0^\infty \cdots \int_0^\infty (t_1^2 + \ldots + t_m^2
+ (t_1 \cdots t_m)^{-2})^{-ns/2} 
\frac{dt_1}{t_1} \cdots \frac{dt_m}{t_m} \ ,
\end{eqnarray}
The argument above starting right after formula~(\ref{eqn:real})
goes through, with the following
changes: take
\begin{eqnarray} \label{eqn:d2}
d(s) = \int_0^\infty \cdots \int_0^\infty (t_1^2 + \ldots + t_m^2
+ (t_1 \cdots t_m)^{-2})^{-ns/2} 
\frac{dt_1}{t_1} \cdots \frac{dt_m}{t_m}
= \frac{\Gamma(s/2)^{m+1}}{2^m (m+1) \Gamma(ns/2)}
\ ,
\end{eqnarray}
replace throughout the term $(t_1^{\delta_1} \cdots t_m^{\delta_m})$
by $(t_1 \cdots t_m)$, 
replace the line containing equation~(\ref{eqn:units}) by 
``Now $1 = |N \epsilon_i| = \prod_{j=1}^{m+1} |\epsilon_i|_j^2$,
and so $|\epsilon_i|_{m+1} = (|\epsilon_i|_1^{\delta_1} \cdots
|\epsilon_i|_m^{\delta_m})^{-1}$'',
and replace the last entry in the $i$-th row of~$M$ by the two
entries
$\Re(\alpha_{i,m+1}) (t_1 \cdots t_m)^{-2}$ and
$\Im(\alpha_{i,m+1}) (t_1 \cdots t_m)^{-2}$.
So equation~(\ref{eqn:seminfinal})
is still valid, but with $d(s)$ given by formula~(\ref{eqn:d2}),
and with the definition of~$\tau(t_1, \ldots, t_m)$ modified
as per the change in~$M$ mentioned above.

Let $V = \int_D\frac{dt_1}{t_1} \cdots \frac{dt_m}{t_m}$.
Then  from equation~(\ref{eqn:seminfinal}) and our 
generalization of Kronecker's limit formula
(formula~(\ref{eqn:kron1})), we get

\begin{thm} \label{thm:zeta}
Recall that $K$ is a number field that is not the rational numbers or
a quadratic imaginary field. With notation as above, in paricular,
taking $d(s)$ to be given by formula~(\ref{eqn:d1})  if $K$ has
a real embedding, and by formula~(\ref{eqn:d2}) if not, we have
\begin{eqnarray*}
w_K  \big(2^{-c} \pi^{-n/2} d_K^{1/2}\big)^s d(s) \Gamma(ns/2)
\zeta_K(s,A)  & =  &
\frac{2V/n}{s-1}
+ (\gamma - \log 4 \pi)V  
\\ & & 
- \frac{2}{n} \int_D
\log \bigg(\prod_{i=1}^{n-1} y_i(t_1, \ldots, t_m)^{i}\bigg) \frac{dt_1}{t_1} \cdots \frac{dt_m}{t_m}
\\ & & 
- 4 \int_D \log g(\tau(t_1, \ldots, t_m))
\frac{dt_1}{t_1} \cdots \frac{dt_m}{t_m} + O(s-1) .
\end{eqnarray*}
\end{thm}
When $n=2$, we recover Hecke's formula~(\ref{eqn:hecke})
and when $n=3$ and $K$ has a complex embedding, we get
Theorem~3 of~\cite{efrat}.
Finally, note that the Dedekind zeta function is the sum of the partial
zeta functions over all ideal classes, so from the formula above,
we get a corresponding formula involving the Dedekind zeta function.

\section{Proof of Theorem~\ref{thm:main}} \label{sec:proof}
We first prove the formula for~$E_n^*(\tau, s)$ and deduce
from it the formula for~$E^\dagger_n(\tau, s)$.
In formula~(\ref{eqn:e*}),
the term corresponding to $m_1=0$  is
\begin{eqnarray} \label{eqn:s1}
S_1 & = &
\pi^{-ns/2} \Gamma(ns/2) 
\sum_{m_2,\ldots, m_n} 
\frac{(\prod_{i=1}^{n-1} y_i^{n-i})^s}
{\parallel (m_2 \ldots m_n) \tau' \parallel^{ns/2}} \nonumber \\
& = &\pi^{-ns/2} \Gamma(ns/2)  \cdot
\bigg(\prod_{i=1}^{n-1} y_i^{(\frac{i}{n-1})}\bigg)^s \cdot
{\sum_{m_2,\ldots, m_n}} 
\frac{(\prod_{i=1}^{n-2} y_i^{n-1-i})^{(\frac{n}{n-1}s)}}
{\parallel (m_2 \ldots m_n) \tau' \parallel^{(n-1)(\frac{n}{n-1}s)/2}} \nonumber \\
& = &  
\bigg(\prod_{i=1}^{n-1} y_i^{(\frac{i}{n-1})}\bigg)^s \cdot
\pi^{-ns/2} \Gamma(ns/2)  
\cdot E_{n-1}\bigg(\tau', \frac{n}{n-1} s\bigg) \nonumber\\
& = &  \bigg(\prod_{i=1}^{n-1} y_i^{(\frac{i}{n-1})}\bigg)^s \cdot
\pi^{-(n-1)(\frac{n}{n-1}s)/2} \Gamma\bigg((n-1)\Big(\frac{n}{n-1}s\Big)/2\bigg)
\cdot E_{n-1}\bigg(\tau', \frac{n}{n-1} s\bigg) \nonumber\\
& = &\bigg(\prod_{i=1}^{n-1} y_i^{(\frac{i}{n-1})}\bigg)^s \cdot 
E^*_{n-1}\bigg(\tau', \frac{n}{n-1} s\bigg) .
\end{eqnarray}

Let
$$S_2 = 
\sum_{m_1 \neq 0} \sum_{m_2, \ldots,  m_n}
\frac{\pi^{-ns/2} \Gamma(ns/2) 
}{\parallel 
(m_1 \ldots m_n) \tau \parallel^{ns/2}} 
,$$
so that 
\begin{eqnarray} \label{eqn:sum}
E_n^*(\tau, s) = S_1 + \bigg(\prod_{i=1}^{n-1} y_i^{n-i}\bigg)^s \cdot S_2.
\end{eqnarray}

Our next goal is to find a suitable expression for~$S_2$, which will
not be achieved till equation~(\ref{eqn:s2''}) below.
We use the formula
\begin{eqnarray} \label{eqn:gamma}
\frac{\pi^{-s} \Gamma(s)}{a^{s}} = \int_0^\infty \exp(- \pi a t) t^{s} \frac{dt}{t}
\end{eqnarray} 
with $a = \parallel 
(m_1 \ldots m_n) \tau \parallel$, and $s$ replaced by~$ns/2$ to get
\begin{eqnarray} \label{eqn:s21}
S_2 = 
\sum_{m_1 \neq 0} \sum_{m_2, \ldots,  m_n} \inti \exp(-\pi t \parallel 
(m_1 \ldots m_n) \tau \parallel) 
t^{ns/2} \frac{dt}{t} \ .
\end{eqnarray} 
For $j = 2, \ldots, n$, let 
\begin{eqnarray} \label{eqn:a}
a_j = (m_1 \tau_{1,1})^2 + \cdots + \bigg(\sum_{i = 1, \ldots, j-1} m_i \tau_{i,j-1}\bigg)^2,
\end{eqnarray}
and recall that 
$b_j = \sum_{i = 1, \ldots, j-1} m_i \tau_{i,j}$,
and $c_j = \tau_{j,j}$.
Then
\begin{eqnarray} \label{eqn:an}
\parallel 
(m_1 \ldots m_n) \tau \parallel & = &
(m_1 \tau_{1,1})^2 + \cdots + 
\bigg(\sum_{i = 1, \ldots, n-1} m_i \tau_{i,n-1}\bigg)^2
+ \bigg( \bigg( \sum_{i = 1, \ldots, n-1} m_i \tau_{i,n} \bigg) + m_n \tau_{n,n}
\bigg)^2 \nonumber \\
& = & a_n + (b_n+ c_n  m_n)^2 , 
\end{eqnarray}
Putting~(\ref{eqn:an}) in~(\ref{eqn:s21}), we get
\begin{eqnarray} \label{eqn:s22}
S_2 & = & 
\sum_{m_1 \neq 0} \sum_{m_2, \ldots,  m_n} \inti \exp(-\pi t a_n) 
\exp(-\pi t(b_n+ c_n m_n)^2) t^{ns/2} \frac{dt}{t} \nonumber \\
& = & 
\sum_{m_1 \neq 0} \sum_{m_2, \ldots,  m_{n-1}} \inti \exp(-\pi t a_n) 
\sum_{m_n}\exp(-\pi t(b_n+ c_n m_n)^2)
t^{ns/2} \frac{dt}{t}
\end{eqnarray}

The Poisson summation formula says that for real numbers~$t, b, c$, with $c \neq 0$,
\begin{eqnarray}\label{eqn:poisson}
 \sum_{m \in \Z} \exp( -\pi t (b + c m)^2) 
 = \frac{1}{c \sqrt{t}} \sum_{m \in \Z} \exp(2 \pi i b m/c) \exp(-\pi m^2/t c^2).
\end{eqnarray}
Using this with $b= b_n$ and $c= c_n$, replacing $m$ by~$m_n$,
and noting that $c_{n}$, being
a diagonal entry of~$\tau$, is always positive, we get
\begin{eqnarray*} \label{eqn:poi1}
\sum_{m_{n}} \exp(-\pi t(b_{n} + c_{n} m_{n})^2)
 =   
\frac{1}{c_{n} \sqrt{t}} \exp(-(\pi t a_{n})) 
\sum_{m_{n}} 
\exp(2 \pi i b_n m_{n}/c_{n}) \exp(-\pi m_{n}^2/t c_{n}^2).
\end{eqnarray*}
Putting this in~(\ref{eqn:s22}), we get 
\begin{eqnarray} \label{eqn:1}
S_2 & =  & 
\frac{1}{c_n} \sum_{m_1 \neq 0} \sum_{m_2, \ldots,  m_{n-1}} \inti \exp(-\pi t a_n) 
\sum_{m_n} \exp(2 \pi i b_n m_n/c_n) \exp(-\pi m_n^2/t c_n^2 )
t^{ns/2 - 1/2} \frac{dt}{t} \nonumber \\
& = & \frac{1}{c_n} 
\sum_{m_1 \neq 0} \sum_{m_2, \ldots,  m_{n}} 
\exp(2 \pi i b_n m_n/c_n)
\inti \exp(-(\pi a_n t   + \pi m_n^2/t c_n^2))
t^{ns/2 - 1/2} \frac{dt}{t}\nonumber \\
& = & \frac{1}{c_n} 
\sum_{m_1 \neq 0, m_n } \exp(2 \pi i b_n m_n/c_n)\nonumber  \\
& & \cdot \inti
\exp(- \pi (m_n^2/c_n^2)/t)
\sum_{m_2, \ldots,  m_{n-2}} 
\sum_{m_{n-1}} \exp(-(\pi a_n t))  
t^{ns/2 - 1/2} \frac{dt}{t}
\end{eqnarray}
Now from equation~(\ref{eqn:a}),
\begin{eqnarray*} 
a_n & = & 
(m_1 \tau_{1,1})^2 + \cdots + 
\bigg(\sum_{i = 1, \ldots, n-2} m_i \tau_{i,n-2}\bigg)^2
+ \bigg( \bigg( \sum_{i = 1, \ldots, n-2} m_i \tau_{i,n-1} \bigg) + m_{n-1} \tau_{n-1,n-1}\bigg)^2 \\
& = & a_{n-1} + (b_{n-1}+ c_{n-1}  m_{n-1})^2 \ .
\end{eqnarray*}
So using formula~(\ref{eqn:poisson}) again as above,
we have
\begin{eqnarray*}
& & \sum_{m_{n-1}} \exp(-(\pi a_n t))  \\
& = & \exp(-(\pi t a_{n-1})) 
\sum_{m_{n-1}} \exp(-\pi t(b_{n-1} + c_{n-1} m_{n-1})^2)\\
& = &  
\frac{1}{c_{n-1} \sqrt{t}} \exp(-(\pi t a_{n-1})) 
\sum_{m_{n-1}} 
\exp(2 \pi i b_n m_{n-1}/c_{n-1}) \exp(-\pi m_{n-1}^2/t c_{n-1}^2).
\end{eqnarray*}

Putting this in~(\ref{eqn:1}),
\begin{eqnarray} \label{eqn:n}
S_2 & =  & \frac{1}{c_{n-1} c_n }
\sum_{m_1 \neq 0, m_n } \exp(2 \pi i b_{n} m_n/c_n)
\inti
\exp(-\pi m_n^2/t c_n^2)
\sum_{m_2, \ldots,  m_{n-2}} 
\frac{1}{\sqrt{t}} \exp(-(\pi t a_{n-1})) \nonumber \\
& & \cdot \sum_{m_{n-1}} 
\exp(2 \pi i b_{n-1} m_{n-1}/c_{n-1}) \exp(-\pi m_{n-1}^2/t c_{n-1}^2)
t^{ns/2 - 2/2} \frac{dt}{t} \nonumber \\
& = & \frac{1}{c_{n-1} c_n }
\sum_{m_1 \neq 0, m_n } \sum_{m_{n-1}} 
\exp(2 \pi i (b_{n} m_n/c_n + b_{n-1} m_{n-1}/c_{n-1}))\nonumber \\
& & \cdot \inti \exp(-\pi (m_n^2/c_n^2 + m_{n-1}^2/c_{n-1}^2)/t)
\sum_{m_2, \ldots,  m_{n-2}}  \exp(-(\pi t a_{n-1})) t^{ns/2 - 2/2} \frac{dt}{t} 
\nonumber \\
& = & \frac{1}{c_{n-1} c_n }
\sum_{m_1 \neq 0, m_n,m_{n-1}} 
\exp(2 \pi i (b_{n} m_n/c_n + b_{n-1} m_{n-1}/c_{n-1})) \nonumber \\
& & \cdot \inti
\exp(-\pi (m_n^2/c_n^2 + m_{n-1}^2/c_{n-1}^2)/t)
\sum_{m_2, \ldots,  m_{n-3}} \sum_{m_{n-2}}
\exp(-(\pi t a_{n-1})) 
t^{ns/2 - 2/2} \frac{dt}{t} 
\end{eqnarray}

Looking at equations~(\ref{eqn:1}) and~(\ref{eqn:n}),
we see that 
repeated use of Poisson summation gives
\begin{eqnarray*}
S_2  =  \frac{1}{\prod_{i=2}^{n} c_i}
\sum_{m_1 \neq 0, m_n, m_{n-1}, \ldots, m_2} 
\exp(2 \pi i d)
\inti
\exp(-\pi m^2/t)
\exp(-(\pi t a_1))
t^{ns/2 - (n-1)/2} \frac{dt}{t} \ ,
\end{eqnarray*}
where recall that $m$ and~$d$ were defined in Theorem~\ref{thm:main}.
Now $a_1 = m_1^2 y^2$ where $y =  \tau_{1,1} = y_1 y_2 \cdots y_{n-1}$. So
\begin{eqnarray*}
\bigg( {\prod_{i=2}^{n} c_i} \bigg) S_2 = 
\sum_{m_1 \neq 0} \sum_{m_2, \ldots,  m_{n}} 
\exp(2 \pi i d) \inti \exp(-( \pi (m_1 y)^2 t +  \pi m^2/t))
t^{n(s - 1)/2 + 1/2} \frac{dt}{t} .
\end{eqnarray*}

Denote the term corresponding $m_2 = \cdots = m_n=0$ by~$S'_2$, i.e., 
$$S'_2= 
\sum_{m_1 \neq 0}
\inti \exp(-\pi (m_1 y)^2 t)
t^{n(s - 1)/2 + 1/2}\frac{dt}{t} \ .$$
Using formula~(\ref{eqn:gamma}), with $s$ replaced by~$n(s - 1)/2 + 1/2$,
we get 
\begin{eqnarray} \label{eqn:s2'}
S'_2 & = & 
\sum_{m_1 \neq 0} 
\frac{\pi^{-(n(s - 1)/2 + 1/2)} \Gamma(n(s - 1)/2 + 1/2)}
{(m_1 y)^{2(n(s - 1)/2 + 1/2)}} \nonumber \\
& = & 
y^{-(n(s - 1) + 1)}
\pi^{-(n(s - 1)/2 + 1/2)} \Gamma(n(s - 1)/2 + 1/2)  \cdot 2
\sum_{m_1 \geq 0} \frac{1}{m_1^{n(s - 1) + 1}}  \nonumber \\
& = & 
2 y^{-(n(s - 1) + 1)}
\pi^{-(n(s - 1)/2 + 1/2)} \Gamma(n(s - 1)/2 + 1/2) \zeta(n(s - 1) + 1)
\end{eqnarray}

Let 
\begin{eqnarray} \label{eqn:s2}
 S''_2 & = & \bigg( {\prod_{i=2}^{n} c_i} \bigg) S_2 - S'_2\\
& = & \sum_{m_1 \neq 0} \sum_{(m_2, \ldots,  m_n) \neq (0,\ldots, 0)} 
\exp(2 \pi i d) \inti \exp(-( \pi (m_1 y)^2 t +  \pi m^2/t))
t^{n(s - 1)/2 + 1/2} \frac{dt}{t} . \nonumber
\end{eqnarray} 
For $a$ and~$b$ positive real numbers, recall the function 
$$K_s(a,b) = \inti \exp(-(a^2 t + b^2/ t)) t^s \frac{dt}{t}$$
Noting that $m \neq 0$ if 
not all $m_2, \ldots, m_n$ are zero,
\begin{eqnarray} \label{eqn:s2''}
S''_2 = 
\sum_{m_1 \neq 0} \sum_{(m_2, \ldots,  m_n) \neq (0,\ldots, 0)} 
\exp(2 \pi i d)
K_{\scriptscriptstyle {n(s - 1)/2 + 1/2}}(\sqrt{\pi} |m_1| y, \sqrt{\pi} |m|)
\end{eqnarray}

From equations~(\ref{eqn:s1}), (\ref{eqn:sum}),  (\ref{eqn:s2'}), (\ref{eqn:s2}), 
and~(\ref{eqn:s2''}), we finally get an expression
for~$E_n^*(\tau, s)$:
\begin{eqnarray}\label{eqn:en}
E_n^*(\tau, s) 
& = & \bigg(\prod_{i=1}^{n-1} y_i^{(\frac{i}{n-1})}\bigg)^s \cdot 
E^*_{n-1}\bigg(\tau', \frac{n}{n-1} s\bigg)  
\nonumber \\
& & 
+ 2 \bigg(\prod_{i=1}^{n-1} y_i^{n-i}\bigg)^s \frac{1}{\prod_{i=2}^{n} c_i}
y^{-(n(s - 1) + 1)}
\pi^{-(n(s - 1)/2 + 1/2)} \Gamma(n(s - 1)/2 + 1/2) \zeta(n(s - 1) + 1)
 \nonumber \\
& & 
+ \bigg(\prod_{i=1}^{n-1} y_i^{n-i}\bigg)^s 
\frac{1}{\prod_{i=2}^{n} c_i}
\sum_{m_1 \neq 0} \sum_{\stackrel{(m_2, \ldots,  m_{n})}{\neq (0,\ldots, 0)}} 
\exp(2 \pi i d)
K_{\scriptscriptstyle {n(s - 1)/2 + 1/2}}(\sqrt{\pi} |m_1| y, \sqrt{\pi} |m|)
\end{eqnarray}

The good thing about the formula above is that it is easy to read off
the polar part and the constant term in each of the summands above, which
is what we do now.
It is known that $K_s$ is an entire function of~$s$, and so all the functions appearing
in the expression
above are holomorphic at~$s=1$ except 
$\zeta(2s-1)$, which has a simple pole at~$s=1$,
and perhaps $E^*_{n-1}\bigg(\tau', \frac{n}{n-1} s\bigg)$.
By induction, $E^*_{n-1}\bigg(\tau', \frac{n}{n-1} s\bigg)$
is also holomorphic except
perhaps when $(\frac{n}{n-1}) s = 1$, and in particular is
homlomorphic at $s=1$. So the first and last
summands on the right side of equation~(\ref{eqn:en}) 
are holomorphic at $s=1$; using the fact
that $K_{1/2}(a,b) = \frac{\sqrt{\pi}}{a} \exp(-2ab)$,
their sum is
\comment{
$$
\bigg(\prod_{i=1}^{n-1} y_i^{n-i-1}\bigg)
\bigg( E^*_{n-1}\bigg(\tau', \frac{n}{n-1}\bigg)  
\sum_{m_1 \neq 0} \sum_{(m_2, \ldots,  m_{n}) \neq (0,\ldots, 0)} 
\exp(2 \pi i d)
\frac{\sqrt{\pi}}{\sqrt{\pi} |m_1| y} \exp(-\sqrt{\pi} |m_1| y \sqrt{\pi} |m|
\bigg)
$$
}
$$
\bigg(\prod_{i=1}^{n-1} y_i^{(\frac{i}{n-1})}\bigg) 
E^*_{n-1}\bigg(\tau', \frac{n}{n-1}\bigg)  + 
\sum_{m_1 \neq 0} \sum_{(m_2, \ldots,  m_{n}) \neq (0,\ldots, 0)} 
\exp(2 \pi i d)
\frac{1}{|m_1|} \exp(- 2 \pi  |m_1| |m| y)
+ O(s-1),
$$
which is $-4 \log g(\tau)+ O(s-1)$.

In order to deal with the second summand 
on the right side of equation~(\ref{eqn:en}),
note that
$$\zeta(n(s - 1) + 1) = \frac{1}{n(s-1)} + \gamma + O(s-1),$$
$$\Gamma(n(s - 1)/2 + 1/2) = \sqrt{\pi} (1 + \frac{n}{2}(\gamma - \log 4)(s-1) + O(s-1)^2,$$
$$\pi^{-(n(s - 1)/2 + 1/2)} 
= \frac{1}{\sqrt{\pi}} ( 1 - \frac{n}{2}\log \pi (s-1) + O(s-1)^2,$$
$$y^{-(n(s - 1) + 1)} 
= y^{-1} ( 1 - n \log y (s-1) + O(s-1)^2, $$
and
$$\bigg(\prod_{i=1}^{n-1} y_i^{n-i}\bigg)^s 
= \bigg(\prod_{i=1}^{n-1} y_i^{n-i}\bigg)^{(s-1) + 1}
= \bigg(\prod_{i=1}^{n-1} y_i^{n-i}\bigg) (1 + \log \bigg(\prod_{i=1}^{n-1} y_i^{n-i}\bigg) (s-1) + O(s-1)^2) .
$$
Using the formulas above, the second summand in 
on the right side of equation~(\ref{eqn:en}) becomes
$$
\frac{2/n}{s-1} + \bigg(\gamma - \log 4\pi 
- \frac{2}{n}\log \bigg(\prod_{i=1}^{n-1} y_i^{i}\bigg)
\bigg)  + O(s-1) 
$$
Using the formulas obtained above
for the three 
summands on the right side of equation~(\ref{eqn:en}),
we get the formula for $E_n^*(\tau, s)$ given in Theorem~\ref{thm:main}.

\bibliographystyle{amsalpha}         

\providecommand{\bysame}{\leavevmode\hbox to3em{\hrulefill}\thinspace}
\providecommand{\MR}{\relax\ifhmode\unskip\space\fi MR }
\providecommand{\MRhref}[2]{%
  \href{http://www.ams.org/mathscinet-getitem?mr=#1}{#2}
}
\providecommand{\href}[2]{#2}

\end{document}